\documentstyle[11pt,leqno]{article}
\def\beq{\begin{equation}}
\def\eeq{\end{equation}}
\newtheorem{th}{Theorem}[section]
\newtheorem{lem}[th]{Lemma}
\newtheorem{pro}[th]{Proposition}

\title{\Large \bf {Quaternionic K\"ahler and hyperK\"ahler manifolds with torsion\\ and
twistor spaces}}
\author{{\sc Stefan Ivanov}
\thanks{Partially supported by Contract MM 809/1998 with the
Ministry of Science and Education of Bulgaria, Contract 586/2002
with the University of Sofia "St. Kl. Ohridski", The European
Contract Human Potential Programme, Research Training Network
HPRN-CT-2000-00101".} \hspace{4mm} and {\sc Ivan Minchev}}
\date{}
\begin{document}
\maketitle
\thispagestyle{empty}
\begin{center}
%Department
%of Mathematics\\ University of Sofia\\ "St. Kl. Ohridski"\\[2mm]
%{\sc The Abdus Salam ICTP Preprint, December, 2001}
\end{center}
\vspace{8mm}

\begin{abstract}
The target space of a (4,0) supersymmetric two-dimensional sigma
model with a Wess-Zumino term has a connection with a totally
skew-symmetric torsion and holonomy contained in Sp(n)Sp(1) (resp.
Sp(n)), QKT (resp. HKT)-spaces.  We study the geometry of QKT and
HKT manifolds and their twistor spaces.  We show that the Swann
bundle of a QKT manifold admits a HKT structure with special
symmetry, if and only if the twistor space of the QKT manifold
admits an almost hermitian structure with totally skew-symmetric
Nijenhuis tensor. In this way we connect two structures arising
from quantum field theories and supersymmetric sigma models with
Wess-Zumino term.
\\[8mm] {\bf Running title:} Quaternionic K\"ahler with
torsion \\[8mm] {\bf Keywords.}  Almost Quaternionic, Hyper Hermitian,
Quaternionic K\"ahler, Torsion, Twistors.
\\[8mm]
\end{abstract}
%\newpage

\section{Introduction and statement of the results}

An almost hyper complex structure on a 4n-dimensional manifold $M$ is
a triple $H=(J_{\alpha}), \alpha=1,2,3$, of almost complex structures
$J_{\alpha}:TM\rightarrow TM$ satisfying the quaternionic identities
$J_{\alpha}^2=-id$ and $J_1J_2=-J_2J_1=J_3$. When each $J_{\alpha}$ is
a complex structure, $H$ is said to be a hyper complex structure on
$M$.

An almost quaternionic structure on $M$ is a rank-3 subbundle $Q
\subset End(TM)$ which is locally spanned by almost hypercomplex
structure $H=(J_{\alpha})$. Such a locally defined triple $H$ will be
called an admissible basis of $Q$. A linear connection $\nabla$ on
$TM$ is called a quaternionic connection if $\nabla$ preserves $Q$,
i.e. $\nabla_X\sigma\in \Gamma (Q)$ for all vector fields $X$ and
smooth sections $\sigma \in \Gamma (Q)$. An almost quaternionic
structure is said to be quaternionic if there is a torsion-free
quaternionic connection. A $Q$-hermitian metric is a Riemannian metric
which is Hermitian with respect to each almost complex structure in
$Q$. An almost quaternionic (resp. quaternionic) manifold with
$Q$-hermitian metric is called an almost quaternionic Hermitian (resp.
quaternionic hermitian) manifold

For $n=1$ an almost quaternionic structure is the same as an
oriented conformal structure and it turns out to be always
quaternionic. When $n\ge 2$, the existence of torsion-free
quaternionic connection is a strong condition which is equivalent
to the 1-integrability of the associated $GL(n,{\bf H})Sp(1)$
structure \cite{Bon,Ob,Sal4}.  If the Levi-Civita connection of a
quaternionic hermitian manifold $(M,g,Q)$ is a quaternionic
connection then $(M,g,Q)$ is called a Quaternionic K\"ahler
manifold (briefly QK manifold). This condition is equivalent to
the statement that the holonomy group of $g$ is contained in
$Sp(n)Sp(1)$ \cite{A1,A2,S1,S2,Ish}.  If on a QK manifold there
exists an admissible basis $(H)$ such that each almost complex
structure $(J_{\alpha})\in(H), \alpha =1,2,3$ is parallel with
respect to the Levi-Civita connection then the manifold is called
hyperK\"ahler (briefly HK). In this case the holonomy group of $g$
is contained in Sp(n).

The various notions of quaternionic manifolds arise in a natural way
from the theory of supersymmetric sigma models as well as in string
theory.  The geometry of the target space of two-dimensional sigma
models with extended supersymmetry is described by the properties of a
metric connection with torsion \cite{HP1,HP2}.  The geometry of (4,0)
supersymmetric two-dimensional sigma models without Wess-Zumino term
(torsion) is a hyperK\"ahler manifold. In the presence of torsion the
geometry of the target space becomes hyperK\"ahler with torsion
(briefly HKT) \cite{HP3}. This means that the complex structures
$J_{\alpha}, \alpha=1,2,3$, are parallel with respect to a metric
quaternionic connection with totally skew-symmetric torsion
\cite{HP3}. Local (4,0) supersymmetry requires that the target space
of two dimensional sigma models with Wess-Zumino term be either HKT or
quaternionic K\"ahler with torsion (briefly QKT)
\cite{Hish} which means that the quaternionic subbundle is parallel
with respect to a metric linear connection with totally skew-symmetric
torsion and the torsion 3-form is of type (1,2)+(2,1) with respect to
all almost complex structures in $Q$.  The target space of
two-dimensional (4,0) supersymmetric sigma models with torsion coupled
to (4,0) supergravity is a QKT manifold \cite{HOP}.

HKT spaces with symmetry (homothety) arise in quantum field
theories. The geometry coming from the Michelson and Strominger's
study of $N=4B$ supersymmetric quantum mechanics with superconformal
$D(2,1,\alpha)$-symmetry is a HKT geometry with a special homothety
\cite{MiS}. These special HKT spaces are studied recently in \cite{PS,PS1}.
It is shown in \cite{PS1} that the special homothety generates an
infinitesimal action of the non-zero quaternions and the quotient
space carries a QKT structure which is of instanton type. Explicitly,
this means that we can find a certain torsion-free quaternionic
connection which induces on the real canonical bundle $k^{{\bf
R}}=\Lambda^{4n}_{{\bf R}}T^*M$ a connection whose curvature is of type
(1,1) with respect to each $J_{\alpha}$. Conversely, for a QKT of
instanton type, one can find (see \cite{PS1}) a HKT structure with
special homothety on the corresponding Swann bundle (a bundle
constructed by A. Swann for QK manifold \cite{Sw}) provided some
nondegeneracy (positivity) conditions are fulfilled.

HKT manifolds are also important in string theory. The number of
surviving supersymmetries in a compactification of a 10-dimensional
string theory on $M$, depends on the number of spinors parallel with
respect to a connection with totally skew-symmetric torsion. This
imposes restrictions on the holonomy group: the spinor representation
of the holonomy group should have a fixed non-trivial spinor. The HKT
geometry is one of the possible models of such a compactification
since the holonomy group of the HKT-connection is a subgroup of
$Sp(n)$. For a more precise discussion concerning parallel spinors and
holonomy of connection with torsion the reader may wish to consult
\cite{Stro,IP,FI,I,FI1}.

The properties of HKT and QKT geometries resemble those of HK and QK
ones, respectively. In particular, HKT \cite{HP3,GP} and QKT
\cite{HOP} manifolds admit twistor constructions with twistor spaces
which have similar properties to those of HK \cite{HKLR} and QK
\cite{S1,S2,S3} assuming some conditions on the torsion
\cite{HP3,HOP,PS1}.  It is shown in \cite{S3,HOP} that the twistor
space of a QKT manifold is always a complex manifold provided that its
dimension is at least 8.  Most of the known examples of QKT manifolds
(e.g. the ones constructed in \cite{homo}) are homogeneous.  However
there is also a large class QKT spaces obtained by conformal
transformations of QK or HK manifolds
\cite{I1}.

The main object of interest in this article are the differential
geometric properties of QKT and HKT manifolds and their
twistor spaces.  We find relations between Riemannian scalar
curvatures of a QKT space which allows us to express sufficient
conditions for a compact 8-dimensional QKT manifold to be QK in terms
of its Riemannian scalar curvatures (Theorem~\ref{thm}).

We consider two almost complex structures $I_1$, $I_2$ on the twistor
space $Z$ over a QKT manifold. The structure $I_{1}$ was originally
constructed in \cite{2}, while $I_{2}$ is constructed in \cite{11} for
the QK case. For QKT the integrability of $I_1$ was established in
\cite{HOP}.

We define a family of Riemannian metrics $h_c$ on $Z$ depending on
a parameter $c$, thus obtaining almost hermitian structures on
$Z$. Investigating the corresponding almost hermitian geometry we
prove that $I_2$ is never integrable and that the Swann bundle of
a QKT manifold admits a HKT structure with special symmetry if and
only if $(Z,h_c,I_2)$ is a $G_1$ manifold according to the
Gray-Hervella classification \cite{14} (Theorem~\ref{th2}).  The
class of $G_1$ manifolds can be viewed as a direct sum of
Hermitian and Nearly K\"ahler manifolds and is characterized by
the requirement that the Nijenhuis tensor should be a 3-form.
These manifolds are of particular interest in physics since they
arise as a target spaces of (2,0)-and (2,2)- supersymmetric sigma
models \cite{P1}. The physical applications also require the
existence of a linear connection $\nabla$ preserving the almost
hermitian structure $(g,J)$ and having a totally skew-symmetric
torsion. The $G_1$-manifolds are precisely the object of interest
since this is the largest class where such a connection exists
\cite{FI}.  Using the integrability of $I_1$ we present new
relations between the different Ricci forms (Theorem~\ref{th1}),
i.e. between the different 2-forms which determine the $Sp(1)$
component of the curvature of the QKT connection.

\

\bigskip

\noindent {\bf Acknowledgements.} The final part of the research
was done during the visit of S.I. at the Abdus Salam International
Centre for Theoretical Physics, Trieste Italy.  S.I. thanks the
Abdus Salam ICTP for providing support and an  excellent research
environment. Both authors thank to J.P.Gauntlett for his very
useful comments.

\section{QKT manifolds}

Let ${\bf H}$ be the quaternions and identify ${\bf H^{n}}={\bf
R^{4n}}$.  To fix notation we assume that ${\bf H}$ acts on ${\bf
H}^{n}$ by right multiplication. This defines an
antihomomorphism $\lambda :\{{\rm unit \:quaternions}\}
\longrightarrow SO(4n),$ where our convention is that $SO(4n)$ acts on
${\bf H}^{n}$ on the left. Denote the image by $Sp(1)$ and let
$I_{0}=\lambda (i), J_{0}=\lambda (j), K_{0}=\lambda (k)$.  The Lie
algebra of $Sp(1)$ is $sp(1)=span\{I_{0},J_{0},K_{0}\}$.  Define
$Sp(n)=\{A\in SO(4n): AB=BA$ for all $B\in Sp(1)\}$.  The Lie algebra
of $Sp(n)$ is $sp(n)=\{A\in so(4n): AB=BA$ for all $B\in sp(1)\}$.
Let $Sp(n)Sp(1)$ be the product of the two groups in $SO(4n)$.
Abstractly, $Sp(n)Sp(1)=(Sp(n)\times Sp(1))/{\bf Z}_{2}$.  The Lie
algebra of the group $Sp(n)Sp(1)$ is isomorphic to $sp(n)\oplus
sp(1)$.

Let $(M,g,(J_{\alpha})\in Q,\alpha=1,2,3)$ be a 4n-dimensional almost
quaternionic manifold with $Q$-hermitian Riemannian metric $g$ and an
admissible basis $(J_{\alpha})$. The K\"ahler form $\Phi_{\alpha}$ of
each $J_{\alpha}$ is defined by $\Phi_{\alpha}=g(\bullet,
J_{\alpha}\bullet )$.  Let
$\nabla$ be a quaternionic connection i.e.
\begin{equation}\label{1}
\nabla J_{\alpha}=-\omega_{\beta}\otimes J_{\gamma} +
\omega_{\gamma}\otimes J_{\beta},
\end{equation}
where the $\omega_{\alpha}, \alpha =1,2,3$ are 1-forms.

Here $(\alpha,\beta,\gamma)$ stands for a cyclic permutation of
$(1,2,3)$.

Let $T(X,Y)=\nabla_XY-\nabla_YX-[X,Y]$ be the torsion tensor of type
(1,2) of $\nabla$. We denote by the same letter the torsion tensor of
type (0,3) given by $T(X,Y,Z)=g(T(X,Y),Z)$.

An almost quaternionic manifold $(M,(H_{\alpha})\in \cal Q)$ is a {\it
QKT manifold} if it admits a hermitian quaternionic structure $(g,
\cal Q)$ and a metric quaternionic connection $\nabla$ (a QKT
connection) with a totally skew symmetric torsion which is a
$(1,2)+(2,1)$-form with respect to each $J_{\alpha},
\alpha=1,2,3$. Explicitly this means that
\begin{equation}\label{a1}
T(X,Y,Z)=T(J_{\alpha}X,J_{\alpha}Y,Z)+T(J_{\alpha}X,Y,J_{\alpha}Z) +
T(X,J_{\alpha}Y,J_{\alpha}Z).
\end{equation}
for all $\alpha=1,2,3$.

It follows that the holonomy group of any QKT-connection is a subgroup
of $Sp(n)Sp(1)$, i.e. the bundle $SO(M)$ of oriented orthonormal frames on
a QKT manifold can be reduced to a principal $Sp(n)Sp(1)$-bundle
$P(M)$ so that the QKT-connection 1-form on $P(M)$ is $sp(n)\oplus
sp(1)$-valued.

Every QKT manifold is a quaternionic manifold \cite{I1}. Poon and
Swann constructed explicitly a quaternionic torsion-free connection
$\nabla^q$ on a QKT in \cite{PS1}. Following \cite{PS1} we will say
that a QKT manifold is of instanton type if the curvature of
$\nabla^q$ on the real canonical bundle is of type (1,1) with respect
to each $J_{\alpha}$.  Conversely, any quaternionic manifold locally
admits a QKT structure \cite{PS1}. Globally it is not necessarily true
that a quaternionic hermitian structure on a quaternionic manifold is
a QKT structure. However, if a QKT structure exists then it is unique
and the torsion 3-form is computed in terms of connection 1-forms
$\omega_{\alpha}$ and the exterior derivative of the K\"ahler forms
\cite{I1}.  One consequence of this is the observation that QKT
structures persist under conformal transformations of the metric
\cite{I1}. For HKT, the local existence of a HKT structure on any
hypercomplex manifold is proved in \cite{GP}. On a HKT manifold the
torsion connection is unique. This fact is a consequence of the
general results in \cite{Ga1} (see also \cite{GP}), which imply that
on a hermitian manifold there exists a unique linear connection with
totally skew-symmetric torsion preserving the metric and the complex
structure. This connection is known as the Bismut connection. Bismut
used this connection \cite{BI} to prove a local index theorem for the
Dolbeault operator on non-K\"ahler manifold.  In the physics
literature, the geometry of this connection is referred to as
KT-geometry. Several non-trivial obstructions to the existence of
(non-trivial) Dolbeault cohomology groups on a compact KT-manifold
were described in \cite{AI,IP}.

Ivanov \cite{I1} introduced the torsion 1-form on a QKT manifold by
the equality
\begin{equation}\label{n1}
t(X)=\frac{1}{2}\sum_{i=1}^{4n}T(J_{\alpha}X,e_i,J_{\alpha}e_i),
\end{equation}
where $\{e_i\}, i=1,\ldots, 4n$ is an orthonormal basis, and showed
that it is independent of $J_{\alpha}$.  It turns out that a QKT
structure is of instanton type if and only if the exterior
differential $dt$ of the torsion 1-form is of type (1,1) with respect
to each $J_{\alpha}$ \cite{PS1}.  We define a {\it balanced QKT
manifold} to be a QKT manifold with a zero torsion 1-form. The first
examples of (compact) balanced HKT manifolds were constructed by Dotti
and Fino \cite{DF}.  On a compact QKT manifold one can show \cite{I1}
that a metric with a coclosed torsion 1-form exists in each conformal
class which supports a QKT structure. Such metrics are called
Gauduchon metrics.

\section{Curvature of QKT manifold}

Let $R=[\nabla,\nabla]-\nabla_{[\;,\;]}$ be the curvature tensor
of type (1,3) of $\nabla$. We denote the curvature tensor of type
(0,4) $R(X,Y,Z,V)=g(R(X,Y)Z,V)$ by the same letter. There are
three Ricci forms and six scalar functions given by $
\rho_{\alpha}(X,Y)=\frac{1}{2}
\sum_{i=1}^{4n}R(X,Y,e_i,J_{\alpha}e_i), {\rm
Scal}_{\alpha,\beta}=
-\sum_{i=1}^{4n}\rho_{\alpha}(e_i,J_{\beta}e_i).$ The Ricci forms
satisfy $n[R(X,Y),J_{\alpha}]=\rho_{\gamma}(X,Y)J_{\beta}
-\rho_{\beta}(X,Y)J_{\gamma}$.  The Ricci tensor $Ric$, the scalar
curvatures $Scal$ and $Scal_{\alpha}$ of $\nabla$ are defined by
$Ric(X,Y)=\sum_{i=1}^{4n}R(e_i,X,Y,e_i),
Scal=\sum_{i=1}^{4n}Ric(e_i,e_i), Scal_{\alpha}= -
\sum_{i=1}^{4n}Ric(e_i,J_{\alpha}e_i)$. We shall denote by $R^g,
Ric^g, \rho_{\alpha}^g, etc.$ the corresponding objects for the
metric $g$, i.e. the same objects taken with respect to the
Levi-Civita connection $\nabla^g$.  We may consider
$(g,J_{\alpha})$ as an almost hermitian structure. Then the tensor
$\rho_{\alpha}^{\ast }(X,Y)=\rho^g_{\alpha}(X,J_{\alpha}Y)$ is
known as the $\ast $-Ricci tensor of the almost hermitian
structure. It is equal to $\rho_{\alpha}^{\ast }(X,Y) = -\sum
_{i=1}^{2n}R^g(e_i,X,J_{\alpha}Y,J_{\alpha}e_i)$ by the Bianchi
identity. The function $Scal^g_{\alpha}$ is known also as the
$*$-scalar curvature. If the $*$-Ricci tensor is a scalar multiple
of the metric then the manifold is said to be $*$-Einstein.  In
general, the $*$-Ricci tensor is not symmetric and the
$*$-Einstein condition is a strong condition. We shall see in the
last section that the $*$-Ricci tensors of a HKT manifold are
always symmetric.

The property of a  QKT structure to be of instanton type can be
expressed in terms of the Ricci forms. Namely, a QKT structure is of
instanton type if and only if
each Ricci form $\rho_{\alpha}$ is of type (1,1) with
respect to $J_{\alpha}$ \cite{PS1}.

We show in this section that the scalar curvature functions are not
independent and define a new scalar invariant, the 'quaternionic
$*$-scalar curvature' of a QKT space. We begin with
\begin{pro}\label{nov}
Let $(M,g,(J_{\alpha})\in \cal Q)$ be a 4n-dimensional QKT
manifold. Then the following identities hold \beq\label{nov1}
\sum_{i=1}^{4n}(\nabla_XT)(J_{\alpha}Y,e_i,J_{\alpha}e_i)=2(\nabla_Xt)Y;
\eeq \beq\label{new1}
\sum_{i,j=1}^{4n}dT(e_j,J_{\alpha}e_j,e_i,J_{\alpha}e_i)= -8\delta t
+8||t||^2-\frac{4}{3}||T||^2, \quad
\sum_{i,j=1}^{4n}dT(e_j,J_{\beta}e_j,e_i,J_{\gamma}e_i)=0.  \eeq
\end{pro}
{\bf Proof.} The formula (\ref{nov1}) follows from (\ref{1}) and the
definition (\ref{n1}) of the torsion 1-form by straightforward
calculations.  To prove (\ref{new1}) we need the following algebraic
\begin{lem}\label{lem1}
For a three form $T$ of type (1,2)+(2,1) with respect to each
$J_{\alpha}$ one has
$$
\sum_{i,j=1}^{4n}g(T(e_i,e_j),T(J_{\gamma}e_i,J_{\beta}e_j))=0, \quad
\sum_{i,j=1}^{4n}g(T(e_i,e_j),T(J_{\beta}e_i,J_{\beta}e_j))=\frac{1}{3}||T||^2
$$
\end{lem}
{\bf Proof of the Lemma.} Put
$A=\sum_{i,j=1}^{4n}g(T(e_i,e_j),T(J_{\gamma}e_i,J_{\beta}e_j))$.  Use
(\ref{a1}) three times to get
\begin{eqnarray}\nonumber
2A&=&
\sum_{i,j,k=1}^{4n}T(e_i,e_j,e_k)(T(J_{\gamma}e_i,J_{\gamma}e_j,J_{\alpha}e_k)-
T(J_{\beta}e_i,J_{\beta}e_j,J_{\alpha}e_k)\\ \nonumber
2A&=&\sum_{i,j,k=1}^{4n}T(e_i,e_j,e_k)(-T(J_{\alpha}e_i,e_j,e_k)+
T(J_{\alpha}e_i,J_{\beta}e_j,J_{\beta}e_k)\\ \nonumber
2A&=&\sum_{i,j,k=1}^{4n}T(e_i,e_j,e_k)(T(J_{\alpha}e_i,e_j,e_k)-
T(J_{\gamma}e_i,J_{\gamma}e_j,J_{\alpha}e_k)\\ \nonumber
\end{eqnarray}
Adding these up yields $6A=0$. The proof of the second identity in the
statement of
Lemma~\ref{lem1} is similar and we omit it.

\

\medskip

We need also the expression of $dT$ in terms of $\nabla$ (see
e.g. \cite{I1,IP,FI}),
\begin{eqnarray}\label{13}
dT(X,Y,Z,U)&=&
{\sigma \atop XYZ}\left\{(\nabla_XT)(Y,Z,U) + g(T(X,Y),T(Z,U)\right\}\\
           &-& (\nabla_UT)(X,Y,Z) + {\sigma \atop XYZ}
\left\{g(T(X,Y),T(Z,U)\right\}, \nonumber
\end{eqnarray}
where ${\sigma \atop XYZ}$ denote the cyclic sum of $X,Y,Z$. Taking
the appropriate trace in (\ref{13}) and applying Lemma~\ref{lem1}, we
obtain the first equality in (\ref{new1}). Finally from (\ref{13})
combined with (\ref{nov1}) and Lemma~\ref{lem1} we get that $
\sum_{i,j=1}^{4n}dT(e_j,J_{\beta}e_j,e_i,J_{\gamma}e_i)=
-4\sum_{i,j=1}^{4n}g(T(e_i,e_j),T(J_{\gamma}e_i,J_{\beta}e_j))=0.$
 \ \hfill {\bf Q.E.D.}
\

\medskip

\begin{pro}\label{orp1}
On a 4n-dimensional $(n>1)$ QKT-manifold we have the equalities:
\beq\label{q1}
Scal_{\alpha,\alpha}=Scal_{\beta,\beta}=Scal_{\gamma,\gamma},\quad
Scal_{\alpha,\beta}=0, \quad Scal_{\alpha}=\frac{1}{2}(dt,\Phi_{\alpha})
\eeq
\end{pro}
{\bf Proof.}  Applying (\ref{nov1}) consequently to (3.45), (3.33) in
\cite{I1}, we obtain \beq\label{nov2}
\frac{2(n-1)}{n}\left\{\rho_{\alpha}(X,J_{\alpha}Y)-
\rho_{\beta}(X,J_{\beta}Y)\right\}
=\sum_{i=1}^{4n}\left\{dT(e_i,J_{\alpha}e_i,X,J_{\alpha}Y) -
dT(e_i,J_{\beta}e_i,X,J_{\beta}Y)\right\}; \eeq
\begin{eqnarray}\label{22}
& &(n-1)\rho_{\alpha}(X,J_{\alpha}Y)=-\frac{n(n-1)}{n+2}Ric(X,Y)
+\frac{n(n-1)}{n+2}(\nabla_Xt)Y\\ & & +
\frac{n}{4(n+2)}\sum_{i=1}^{4n}\left\{(n+1)dT(X,J_{\alpha}Y,e_i,J_{\alpha}e_i)-
dT(X,J_{\beta}Y,e_i,J_{\beta}e_i)
-dT(X,J_{\gamma}Y,e_i,J_{\gamma}e_i)\right\}.\nonumber
\end{eqnarray}
Now take the appropriate trace in (\ref{nov2}), and use (\ref{new1})
to get $Scal_{\alpha,\alpha}=Scal_{\beta,\beta},\quad
Scal_{\alpha,\beta}=0$.  The last equality in (\ref{q1}) is a direct
consequence of (\ref{nov1}), $Scal_{\alpha,\beta}=0$ and
(\ref{22}). \hfill {\bf Q.E.D.}

\

\medskip

\noindent
{\bf Definition.} {\em The three coinciding traces of the Ricci forms on a
$4n$ dimensional QKT manifold $(n>1)$, give a well-defined global
function. We call this function the {\bf quaternionic scalar curvature
of the QKT connection} and denote it by $Scal_Q:=
Scal_{\alpha,\alpha}$.}

\

\medskip

\begin{pro}\label{orp2}
On a 4n-dimensional $(n>1)$ QKT manifold we have \beq\label{pq1}
Scal^g_{\alpha}=Scal^g_{\beta}=Scal^g_{\gamma}= Scal_Q - \delta t +
||t||^2 -\frac{1}{12}||T||^2, \quad
Scal^g_{\alpha,\beta}=Scal_{\gamma}=\frac{1}{2}(dt,\Phi_{\gamma}).
\eeq
\end{pro}
{\bf Proof.} We follow \cite{I1,IP}. The curvature $R^g$ of the
 Levi-Civita connection is related to $R$ via
\begin{eqnarray}\label{15}
R^g(X,Y,Z,U) &=& R(X,Y,Z,U) - \frac{1}{2} (\nabla_XT)(Y,Z,U)
+\frac{1}{2} (\nabla_YT)(X,Z,U)\nonumber \\
             &-&
\frac{1}{2}g(T(X,Y),T(Z,U)) - \frac{1}{4}g(T(Y,Z),T(X,U)) -
\frac{1}{4}g(T(Z,X),T(Y,U)).
\end{eqnarray}
Taking the  traces in (\ref{15}) and using (\ref{n1}) we obtain
\begin{eqnarray}\label{rn2}
\rho^g_{\alpha}(X,J_{\alpha}Y)&=&\rho_{\alpha}(X,J_{\alpha}Y)-\frac{1}{2}(\nabla_Xt)Y
-\frac{1}{2}(\nabla_{J_{\alpha}Y}t)J_{\alpha}X\\ \nonumber
&+&\frac{1}{2}t(J_{\alpha}T(X,J_{\alpha}Y))
+\frac{1}{4}\sum _{i=1}^{4n} g\left(T(X,e_i),T(J_{\alpha}Y,J_{\alpha}e_i)\right),
\end{eqnarray}
To finish take the appropriate traces in (\ref{rn2}) and apply Lemma~\ref{lem1}
and Proposition~\ref{orp1}. \ \hfill {\bf Q.E.D.}

\

\medskip

\noindent
{\bf Definition.}  {\em The three coinciding traces of the Riemannian
Ricci forms on a $4n$ dimensional QKT manifold $(n>1)$, give a
well-defined global function. We call this function the {\bf
quaternionic $*$-scalar curvature} and denote it by $Scal^g_Q :=
Scal^g_{\alpha}$.}

\

\medskip

\begin{pro}\label{orp4}
On a 4n-dimensional ($n>1$) QKT manifold $(M,g,\cal Q)$ the scalar
curvatures are related by
\begin{eqnarray}
\nonumber Scal^g &=& \frac{n+2}{n}Scal_Q -3\delta t +2||t||^2
-\frac{1}{12}||T||^2,\\ \nonumber Scal^g_Q &=& Scal_Q - \delta t
+||t||^2 -\frac{1}{12}||T||^2,\\ \nonumber Scal &=&
\frac{n+2}{n}Scal_Q -3\delta t +2||t||^2 -\frac{1}{3}||T||^2.\nonumber
\end{eqnarray}
\end{pro}
{\bf Proof.} We derive from (\ref{15}) that
\begin{eqnarray}\label{r5}
Ric^g(X,Y) &=& Ric(X,Y) + \frac{1}{2}\delta T(X,Y) + \frac{1}{4}\sum
_{i=1}^{2n} g\left(T(X,e_i),T(Y,e_i)\right), \\ & & Scal^g=Scal +
\frac{1}{4}||T||^2.\nonumber
\end{eqnarray}
Take the trace in (\ref{22}) and use Lemma~\ref{lem1} to get the first
equality of the proposition. The second equality is already proved in
Proposition~\ref{pq1}. The last one is a consequence of (\ref{r5}) and
the already proven first equality in the proposition.  \hfill {\bf
Q.E.D.}

\

\medskip

\noindent
As a consequence of the above result, we get

\begin{th}\label{th11}
Let $(M,g,\cal Q)$ be a compact 4n-dimensional $(n>1)$ QKT
manifold. Then \beq\label{ma1} \int_M(Scal^g-Scal^g_Q-
\frac{2}{n}Scal_Q)\,dV\ge 0.  \eeq The equality in (\ref{ma1}) is
attained if and only if the QKT structure is balanced.
\beq\label{ma2} \int_M(Scal^g-2Scal^g_Q- \frac{2-n}{n}Scal_Q)\,dV\ge
0.  \eeq The equality in (\ref{ma2}) is attained if and only if the
QKT structure is quaternionic K\"ahler.
\end{th}
{\bf Proof.} Proposition~\ref{orp4} implies \hspace{4mm}
$Scal^g-Scal^g_Q- \frac{2}{n}Scal_Q=-2\delta t +2||t||^2$,
\beq\label{qaw} Scal^g-2Scal^g_Q- \frac{2-n}{n}Scal_Q=-\delta t
+\frac{1}{12}||T||^2.  \eeq Integrating the last two equalities over
$M$ we get the proof.\hfill {\bf Q.E.D.}

\

\medskip

\noindent
{\bf Remark 1.} From the proof of Theorem~\ref{th11}, it is clear that
the statement of the theorem will still hold for a non-compact QKT,
provided that $t$ is a coclosed 1-form.

\

\medskip

\noindent
Applying Theorem~\ref{th11} to an 8-dimensional QKT manifold, i.e. take
$n=2$ in (\ref{qaw}), we get the main result of this section

\begin{th}\label{thm}
Let $(M,g,\cal Q)$ be an 8-dimensional compact connected QKT
manifold. Then
\begin{itemize}
\item[a)] $(M,g,\cal Q)$ is a quaternionic K\"ahler manifold if and
only if
$$\int_M(Scal^g-2Scal^g_Q)\,dV= 0.$$
\item[b)] $(M,g,\cal Q)$ is a locally hyperK\"ahler manifold if and
only if the Riemannian scalar curvature and the quaternionic
$*$-scalar curvature both vanish.
\end{itemize}
In particular, any compact 8-dimensional QKT manifold with a flat metric
is flat locally hyperK\"ahler and therefore is
covered by a hyperK\"ahler torus.
\end{th}

We finish this section with the following
\begin{th}\label{20}
A $4n$-dimensional QKT manifold is of instanton type if and only if each $*$-Ricci tensor is symmetric.
\end{th}
{\bf Proof.} First, we observe that on a QKT manifold $(M,g,\cal Q)$
the (2,0)+(0,2)-parts of $\rho^g_{\alpha}, \rho_{\alpha}, dt$ with
respect to $J_{\alpha}$ are related by the equality \beq\label{sup}
\rho^g_{\alpha}(X,J_{\alpha}Y)+\rho^g_{\alpha}(J_{\alpha}X,Y)=
\rho_{\alpha}(X,J_{\alpha}Y)+\rho_{\alpha}(J_{\alpha}X,Y) -\frac{1}{2}
\left(dt(X,Y)-dt(J_{\alpha}X,J_{\alpha}Y)\right).  \eeq Indeed, put
$B(X,Y)=\sum _{i=1}^{4n}
g\left(T(X,e_i),T(J_{\alpha}Y,J_{\alpha}e_i)\right)$.  The tensor $B$
is symmetric since the (1,2)+(2,1)-type property of $T$ leads to the
expression\\ $2B(X,Y)=\sum _{i,j=1}^{4n}
\left(T(X,e_i,e_j),T(Y,e_i,e_j)-
T(X,e_i,e_j),T(Y,J_{\alpha}e_i,J_{\alpha}e_j)\right)$ which is clearly
symmetric.  Then the skew-symmetric part of (\ref{rn2}) gives
(\ref{sup}), where we used (\ref{a1}) and the equality
$d^{\nabla}t(X,Y):=(\nabla_Xt)Y-(\nabla_Yt)X=dt(X,Y)-t(T(X,Y))$.
Computations in \cite{PS1} show the identity
$\rho_{\alpha}(X,J_{\alpha}Y)+\rho_{\alpha}(J_{\alpha}X,Y) =
-\frac{n}{2}
\left(dt(X,Y)-dt(J_{\alpha}X,J_{\alpha}Y)\right)$. Consequently,
(\ref{sup}) gives
$\rho^g_{\alpha}(X,J_{\alpha}Y)+\rho^g_{\alpha}(J_{\alpha}X,Y) =
-\frac{n+1}{2} \left(dt(X,Y)-dt(J_{\alpha}X,J_{\alpha}Y)\right)$.
Hence, $M$ is of instanton type if and only if $\rho^g_{\alpha}$ is of
type (1,1) with respect to $J_{\alpha}$. The latter property is
equivalent to the condition that the corresponding $*$-Ricci tensor
$\rho_{\alpha}^*$ is symmetric. \hfill {\bf Q.E.D.}

\section{Twistor space of QKT manifolds}

In this section we adapt the setup from \cite{19,7} to incorporate a
totally skew-symmetric torsion. Our discussion is very close to that
of \cite{AGI}.

Let $(M,g)$ be a $4n$-dimensional QKT
manifold and let $\pi :P(M)\longrightarrow M$ be the natural
projection.  For each $u\in P(M)$ we consider the linear
isomorphism $j(u)$ on $T_{\pi (u)}M$ defined by
$j(u)=uJ_{0}u^{-1}$.
It is easy to see that $j(u)^{2}=-id$ and
$g(j(u)X,j(u)Y)=g(X,Y)$ for all $X,Y \in T_{\pi (u)}M$, i.e.
$j(u)$ is an orthogonal complex structure at $\pi (u)$.  For
each point $p\in M$ we define
$Z_{p}(M)=\{j(u):  u\in P(M), \pi (u)=p\}$.  In other words,
$Z_{p}(M)$ is the space of all orthogonal complex structures
in the tangent space $T_{p}M$ which are compatible with the
QKT structure.

We put $Z=\bigcup_{p\in M}Z_{p}(M)$.
Let $H=Sp(n)Sp(1)\bigcap U(2n)$.  There is a bijective
correspondence between the symmetric space
$Sp(n)Sp(1)/H=Sp(1)/U(1)={\bf CP^1}={\bf S^{2}}$ and
$Z_{p}(M)$ for every $p \in M$.  So we can consider $Z$ as
the associated fibre bundle of $P(M)$ with standard fibre
$Sp(n)Sp(1)/H={\bf CP^1}$.  Hence, $P(M)$ is a principal
fibre bundle over $Z$ with structure group $H$ and
projection $j$.  If $\pi _{1}:Z\longrightarrow M$ is the
projection, we have that $\pi _{1}\circ j=\pi $.
We consider the symmetric space $Sp(n)Sp(1)/H$.  We have
the following Cartan decomposition $sp(n)\oplus sp(1)=h\oplus m,$ where
$h=\{A\in sp(n)\oplus sp(1):  AJ_{0}=J_{0}A\}=(sp(n)\oplus
sp(1))\bigcap u(2n)$ is the Lie algebra of $H$ and
$m=\{A\in sp(n)\oplus sp(1):  AJ_{0}=-J_{0}A\}$.
It is clear that $m$ is generated by $I_{0}$ and $K_{0}$,
i.e. $m=span\{I_{0},K_{0}\}.$
Hence, if $A\in m$ then $J_{0}A\in m$.
Let $( , )$ be the inner product in $gl(4n,{\bf R})$ defined
by
$(A,B)=trace(AB^{t})=\sum_{i=1}^{4n}<Ae_{i},Be_{i}> {\rm
for} A,B\in gl(4n,{\bf R})$, where $< , >$ is  the canonical
inner product in ${\bf R^{4n}}$.
It is clear that $sp(n)\bot sp(1)$ and $I_{0},J_{0},K_{0}$
form an orthogonal basis of $sp(1)$ with
$(I_{0},I_{0})=(J_{0},J_{0})=(K_{0},K_{0})=4n$. Hence,
$h\bot m$.

Let $u\in P(M)$ and $Q_{u}$ is the horizontal subspace of
the tangent space $T_{u}P(M)$ induced by the QKT-
connection on M (\cite{15} ). The vertical space is
$h^{\ast }_{u}\oplus m^{\ast }_{u}$, where
$h^{\ast }_{u}=\{A^{\ast }_{u}: A\in h\},
m^{\ast }_{u}=\{A^{\ast }_{u}: A\in m\}$. Hence,
$T_{u}P(M)=h^{\ast }_{u}\oplus m^{\ast }_{u}\oplus Q_{u}$.
For each $u\in P(M)$ we put
$V_{j(u)}=j_{\ast u}(h^{\ast }_{u}\oplus m^{\ast }_{u}),
H_{j(u)}=j_{\ast u}Q_{u}$. Thus we obtain the vertical and
horizontal distributions $V$ and $H$ on $Z$. Since $P(M)$ is
a principal fibre bundle over $Z$ with structure group $H$ we
have $Ker j_{\ast u}=h^{\ast }_{u}$. Hence
$V_{j(u)}=j_{\ast u}m^{\ast }_{u}$ and
$j_{\ast u|m^{\ast }_{u}\oplus Q_{u}}:m^{\ast }_{u}\oplus
Q_{u}\longrightarrow T_{j(u)}Z$ is an isomorphism.
We define almost complex structures $I_{1}$ and $I_{2}$ on
$Z$ by
\begin{eqnarray}\label{2.2}
& &I_{1}j_{\ast u}A^{\ast }=j_{\ast u}(J_{0}A)^{\ast }, \qquad
I_{2}j_{\ast u}A^{\ast }=-j_{\ast u}(J_{0}A)^{\ast }\\
& &I_ij_{\ast u}B(\xi )=j_{\ast u}B(J_{0}\xi ), \qquad i=1,2,\nonumber
\end{eqnarray}
for $u \in P(M), A \in m, \xi \in {\bf R^{4n}}.$ For twistor bundles
of 4-manifolds the almost complex structure $I_{1}$ is introduced in
\cite{2} and the almost complex structure $I_{2}$ is introduced in
\cite{11} in terms of the horizontal spaces of the Levi-Civita
connection. The almost complex structure $I_1$ for QK, HKT and QKT
manifolds was constructed in \cite{S1,HP3,HOP}, respectively, where it
is proved that it is actually integrable.  For every $c > 0$ a
Riemannian metric $h_{c}$ on Z is defined by
\begin{eqnarray}\label{2.5}
& &h_{c}(j_{\ast u}A^{\ast },j_{\ast u}B^{\ast })=c^{2}(A,B), \qquad
h_{c}(j_{\ast u}A^{\ast },j_{\ast u}B(\xi ))=0\\
& &
h_{c}(j_{\ast u}B(\xi ),j_{\ast u}B(\eta ))=<\xi ,\eta >,\nonumber
\end{eqnarray}
for $u\in P(M), A,B\in m, \xi ,\eta \in {\bf R^{4n}}.$\\
It is clear that $(I_i,h_{c}), i=1,2$, determine two
families of almost hermitian structures on $Z$.

In the QK case the properties of the almost hermitian geometry of
$(I_i,h_{c}), i=1,2$ are considered in \cite{DM,AGI}. Below we follow
\cite{AGI} making the necessary modifications required by the presense
of torsion.

We split the curvature of a QKT-connection into $sp(n)$-valued part
$R'$ and $sp(1)$-valued part $R''$ following the classical scheme (see
e.g. \cite{AM,Ish,Bes})
\begin{pro}\label{p1}
The curvature of a QKT manifold splits as follows
\begin{eqnarray}\nonumber
R(X,Y)&=&R'(X,Y)+{1\over
  2n}(\rho_1(X,Y)J_1+\rho_2(X,Y)J_2+\rho_3(X,Y)J_3),\\ \nonumber & &
  [R'(X,Y),J_\alpha]=0,\ \ \ \alpha=1,2,3.
\end{eqnarray}
\end{pro}

We denote by $A^{\ast }$
(resp.  $B(\xi )$) the fundamental vector field (resp.  the
standard horizontal vector field) on $P(M)$ corresponding to
$A\in sp(n)\oplus sp(1)$ (resp.  $\xi \in {\bf R^{4n}}$).
Let $\Omega, \Theta$ be the curvature 2-form and the torsion 2-form for the
QKT-connection on $P(M)$, respectively (\cite{15} ).

We  shall  denote  the   splitting   of   the
$sp(n)\oplus  sp(1)$-valued  curvature  2-form  $\Omega$  on
$P(M)$, corresponding to Proposition~\ref{p1},
by $\Omega =\Omega'  +\Omega''$, where $\Omega' $ is a $sp(n)$-valued
 2-form and  $\Omega''$  is  a
$sp(1)$-valued form.  Explicitly, we have
$
\Omega ''= \Omega ''_{1}I_{0} + \Omega''_{2}J_{0} + \Omega''_{3}K_{0},
$
where $\Omega''_{\alpha}, \alpha=1,2,3$, are 2-forms.  If
$\xi, \eta, \zeta \in {\bf R^{4n}}$, then the 2-forms
$\Omega''_{\alpha}, \alpha=1,2,3$, are given by
\beq\label{e1}
\Omega''_{\alpha}(B(\xi ),B(\eta ))=\frac{1}{2n}\rho_{\alpha}(X,Y), \quad X=u(\xi), Y=u(\eta).
\eeq
Since $T$ is 3-form of type (1,2)+(2,1), the torsion 2-form $\Theta$ has the properties
\begin{eqnarray}\label{t1}
& & <\Theta_u(B(\xi),B(\eta)),\zeta>=-<\Theta_u(B(\xi),B(\zeta)),\eta>,\nonumber \\
& & <\Theta_u(B(\xi),B(\eta)),\zeta>=<\Theta_u(B(J_0\xi),B(J_0\eta)),\zeta>\\
& &+<\Theta_u(B(J_0\xi),B(\eta)),J_0\zeta>+<\Theta_u(B(\xi),B(J_0\eta)),J_0\zeta> \nonumber
\end{eqnarray}

Let $F_{i}(X,Y,Z) = h_{c}((D_{X}I_i)Y,Z), i=1,2$,
where $D$ is the covariant derivative of the Levi-Civita
connection of $h_{c}$. We denote by K the curvature tensor
of $h_{c}$.

In the rest of the paper
 $A,B,C,D{\in m}, \ \xi,\eta,\zeta,\tau\in\textbf{R}^{4n}$.

The calculations made  in  \cite{AGI}  for  the
twistor space over QK manifold can be  performed  in
our case by taking into account the torsion and their properties. In this way,
using (\ref{2.2}), (\ref{2.5}) and (\ref{t1}), we obtain our technical tools, namely
\begin{pro} The  next equalities
 hold at  $u\in P(M)$:
\begin{eqnarray}\label{2.7}
& &F_i(j_{\ast u}A^{\ast },j_{\ast u}B^{\ast },j_{\ast u}C^{\ast
})=0, \quad F_i(j_{\ast u}A^{\ast },j_{\ast u}B^{\ast },j_{\ast
u}B(\xi ))=0, \quad i=1,2,\nonumber \\ & & F_i(j_{\ast u}A^{\ast
},j_{\ast u}B(\xi ),j_{\ast u}B(\eta )) = \frac{c^{2}}{2}(A,\Omega
(B(J_{0}\xi ),B(\eta )) +\nonumber \\ & & \hspace{50mm}+
\frac{c^{2}}{2}(A,\Omega (B(\xi ),B(J_{0}\eta )) +  2<AJ_{0}\xi
,\eta >, \quad i=1,2,\nonumber \\ & & F_i(j_{\ast u}B(\xi
),j_{\ast u}A^{\ast },j_{\ast u}B^{\ast })=0, \quad i=1,2,  \\
\nonumber & & F_{2}(j_{\ast u}B(\xi ),j_{\ast u}A^{\ast },j_{\ast
u}B(\eta )) = \frac{c^{2}}{2}(A,J_{0}\Omega (B(\xi ),B(\eta ))
+\frac{c^{2}}{2}(A,\Omega (B(\xi ),B(J_{0}\eta )),\\ \nonumber & &
F_{1}(j_{\ast u}B(\xi ),j_{\ast u}A^{\ast },j_{\ast u}B(\eta )) =
 \frac{c^{2}}{2}(J_{0}A,\Omega (B(\xi ),B(\eta ))
+\frac{c^{2}}{2}(A,\Omega (B(\xi ),B(J_{0}\eta )),\\ \nonumber
& &
F_i(j_{\ast u}B(\xi ),j_{\ast u}B(\eta ),j_{\ast u}B(\zeta
))=-\frac{1}{2}<\Theta(B(\xi),B(J_0\eta)),\zeta>-\\ \nonumber
& & \hspace{50mm} -
\frac{1}{2}<\Theta(B(\xi),B(\eta)),J_0\zeta>, \quad i=1,2.
\end{eqnarray}
\end{pro}
For the curvature tensor $K$ of $h_c$ we have
\begin{pro}
The following equalities hold
at any  $u\in P(M)$:
\begin{eqnarray}\label{2.9}
& & K(j_{\ast u}A^{\ast },j_{\ast u}B^{\ast },j_{\ast u}C^{\ast
},j_{\ast u}D^{\ast })=-c^{2}([A,B],[C,D]), \nonumber \\
& & K(j_{\ast u}A^{\ast },j_{\ast u}B^{\ast },j_{\ast u}C^{\ast
},j_{\ast u}B(\xi ))=0, \nonumber \\
& & K(j_{\ast u}A^{\ast },j_{\ast u}B(\xi ),j_{\ast u}B^{\ast
},j_{\ast u}B(\eta ))=\frac{c^{2}}{2}([A,B],\Omega (B(\xi ),
B(\eta))-\nonumber \\
& & \qquad -\frac{c^{4}}{4}(B,\Omega (B(\xi
),B(e_{i}))(u))(A,\Omega (B(\eta
),B(e_{i})),\nonumber \\
& &
K(j_{*u}B(\xi),j_{*u}B(\eta),j_{*u}B(\zeta),j_{*u}A^*)={c^2\over 2}
(A,B(\zeta)\Omega(B(\xi),B(\eta)))\nonumber \\
& & \qquad
-{c^2\over 4}(A,\Omega([B(\eta),B(\zeta)],B(\xi)))
-{c^2\over 4}(A,\Omega([B(\zeta),B(\xi)],B(\eta))),  \\
& &
K(j_{*u}B(\xi),j_{*u}B(\eta),j_{*u}B(\zeta),j_{*u}B(\tau))=
<\Omega(B(\xi),B(\eta))\zeta,\tau> \nonumber \\
& & \qquad
-{c^2\over 4}(\Omega_{m}(B(\xi),B(\tau)),\Omega_{m}(B(\eta),B(\zeta)))
+{c^2\over 4}(\Omega_{m}(B(\xi),B(\zeta)),\Omega_{m}(B(\eta),B(\tau))) \nonumber \\
& & \qquad
+{c^2\over 2}(\Omega_{m}(B(\xi),B(\eta)),\Omega_{m}(B(\zeta),B(\tau)))
-{1\over 4}<\Theta(B(\xi),B(\tau)),\Theta(B(\eta),B(\zeta))> \nonumber \\
& &\qquad
+{1\over 4}<\Theta(B(\xi),B(\zeta)),\Theta(B(\eta),B(\tau))>
-{1\over 2}<\Theta(B(\xi),B(\eta)),\Theta(B(\zeta),B(\tau))> \nonumber \\
& &\qquad
-{1\over 2}<B(\xi)\Theta(B(\eta),B(\zeta)),\tau>
+{1\over 2}<B(\eta)\Theta(B(\xi),B(\zeta)),\tau>, \nonumber
\end{eqnarray}
\end{pro}
 where $\Omega _{m}$ denotes the $m$-component
of $\Omega $.

\section{Almost Hermitian geometry of $(Z,h_c,I_i)$}

Let $(M,g,J)$ be a $2n$-dimensional almost Hermitian
manifold and let $F(X,Y,Z)=g((\nabla^g_{X}J)Y,Z),$ where
$\nabla^g $ is the covariant differentiation of the
Levi-Civita connection on $M$.
We recall the definition of some classes
according to the
Gray-Hervella classification \cite{14} in terms of the
notations we use:
$(M,g,J)$ is K\"ahler if $F=0$; \hspace{1mm} Hermitian if $H(X,Y,Z)=F(X,Y,Z)-F(JX,JY,Z)=0$;
\hspace{1mm} semi-K\"ahler if $tr F =0$, \hspace{1mm} quasi-K\"ahler if $F(X,Y,Z)+F(JX,JY,Z)=0$;
\hspace{1mm} nearly K\"ahler if $F(X,Y,Z)+F(Y,X,Z)=0$; \hspace{1mm}  almost K\"ahler if
$F(X,Y,Z)+F(Y,Z,X)+F(Z,X,Y)=0$; $G_1$ manifold if $\psi(X,Y,Z)=F(X,Y,Z)+F(Y,X,Z)-F(JX,JY,Z)-F(JY,JX,Z)=0.$
\begin{th}\label{th1}
On a QKT mamnifold $(M^{4n},g,(J_{\alpha})\in \cal Q)$ the (2,0)+(0,2)-parts of the Ricci forms
$\rho_{\alpha}$, $\rho_{\beta}$ with respect to $J_{\gamma}$ coincide.
\end{th}
{\bf Proof.} We claim the following identities hold
\beq\label{ric1}
\rho_{\alpha}(J_{\beta}X,J_{\beta}Y)-\rho_{\alpha}(X,Y)=
\rho_{\gamma}(J_{\beta}X,Y)+\rho_{\gamma}(X,J_{\beta}Y).
\eeq
Indeed, consider the almost hermitian structure $(h_c,I_1)$ on the twistor space
$Z$. The almost complex structure $I_1$ is integrable \cite{HOP}. Therefore
$(Z,h_c,I_1)$ is a Hermitian manifold. We calculate taking into account
(\ref{2.7}) that
\begin{eqnarray}\nonumber
& &0=\frac{2}{c^2}\left(F_1(j_{*u}A^*,j_{*u}B(\xi),j_{*u}B(\eta))-
F_1(I_1(j_{*u}A^*),I_1(j_{*u}B(\xi)),j_{*u}B(\eta))\right)=\\ \nonumber
& &(J_0A,\Omega(B(\xi),B(\eta))
+(A,\Omega(B(J_0\xi),B(\eta))+(A,\Omega(B(\xi),B(J_0\eta))
-(J_0A,\Omega(B(J_0\xi),B(J_0\eta)).
\end{eqnarray}

The $m$ component $\Omega_m$ is given by
$2n\Omega_m(B(\xi),B(\eta))=\rho_1(B(\xi),B(\eta))I_0
+\rho_3(B(\xi),B(\eta))K_0$. The last two equalities imply
$$\rho_1(B(J_0\xi),B(\eta))+\rho_1(B(\xi),B(J_0\eta))=
\rho_3(B(J_0\xi),B(J_0\eta))-\rho_3(B(\xi),B(\eta))$$
The proof is completed by puting $J_1=I_0, J_2=J_0, J_3=K_0$.  \hfill {\bf Q.E.D.}

\

\medskip

\noindent
We recall the notion of a Swann bundle \cite{Sw,PPS,PS1}. On a 4n
dimensional QKT manifold $M$ it is defined by ${\cal U}(M) =
(P\times_{Sp(n)Sp(1)}{\bf H}^*)/\{\pm\}$, where ${\bf H}^*$ are the
nonzero quaternions. It carries a hypercomplex structure
\cite{Sw,PPS,PS1}. If $M$ is of instanton type and the condition
\begin{eqnarray}\label{ric3}
\rho_{\alpha}(J_{\alpha}X,Y)+\rho_{\alpha}(J_{\gamma}X,J_{\beta}Y)
=\frac{1}{c^2}g(X,Y)
\end{eqnarray}
holds then ${\cal U}(M)$ has a HKT structure with special homothety
\cite{PS1}. This generalizes the Swann result stating that if $M$ is
QK then ${\cal U}(M)$ carries a HK structure \cite{Sw}. An alternative
construction \cite{HOP} of a HK structure on the Swann bundle over QKT
utilizes the assumption that $dT$ is a (2,2)-form with
respect to each $J_{\alpha}$.
\begin{th}\label{th2}
Let $(M^{4n},g,(J_{\alpha})\in \cal Q)$ be a QKT manifold with twistor
space $Z$.
\begin{itemize}
\item[a)] The almost complex structure $I_2$ on $Z$ is never integrable.
\item[b)] The space $(Z,h_c,I_2)$ is a $G_1$ manifold if and only if
the Swann bundle admits a HKT structure with special homothety.
\item[c)] The spaces $(Z,h_c,I_i), \quad i=1,2$ are semi-K\"ahler
manifolds if and only if the QKT structure is balanced, i.e. if $t=0$.
\item[d)] If $(Z,h_c,I_2)$ is quasiK\"ahler, almost K\"ahler, nearly
K\"ahler or $(Z,h_c,I_1)$ is K\"ahler then the torsion is zero and $M$
is a QK manifold.
\end{itemize}
\end{th}
{\bf Proof.} The almost complex structure $I_2$
is integrable if and only if $H=0$. We obtain using (\ref{2.7}) that
\begin{eqnarray}\nonumber
& &H(j_{*u}B(\xi),j_{*u}A^*,j_{*u}B(\eta))={c^2\over 2}(J_0A,\Omega(B(J_0\xi),B(J_0\eta))) -
{c^2\over 2}(J_0A,\Omega(B(\xi),B(\eta)))\\ \nonumber
& &\hspace{50mm} +{c^2\over 2}(A,\Omega(B(J_0\xi),B(\eta)))
+{c^2\over 2}A,\Omega(B(\xi),B(J_0\eta))) \\ \nonumber
& &H(j_{*u}A^*,j_{*u}B(\xi),j_{*u}B(\eta)) =4<AJ_0\xi,\eta> +
H(j_{*u}B(\xi),j_{*u}A^*,j_{*u}B(\eta)) \nonumber
\end{eqnarray}
 Hence, $<AJ_0\xi,\eta>=0$ which is impossible.

For b),  $(Z,h_c,I_2)$ is a $G_1$-manifold if and only if
$\psi=0$. We get by (\ref{2.7}) that the only non-zero term of $\psi$ is
\begin{eqnarray}\nonumber
& &\psi(j_{*u}A^*,j_{*u}B(\xi),j_{*u}B(\eta))=
4<AJ_0\xi,\eta>-c^2(J_0A,\Omega(B(\xi),B(\eta)))+\\ \nonumber
& &c^2(A,\Omega(B(J_0\xi),B(\eta)))
+c^2(A,\Omega(B(\xi),B(J_0\eta)))+c^2(J_0A,\Omega(B(J_0\xi),B(J_0\eta))) \nonumber
\end{eqnarray}
 The last equality is equivalent to
\begin{eqnarray}\label{ric2}
c^2\rho_1(B(J_0\xi),B(\eta))+c^2\rho_1(B(\xi),B(J_0\eta))&=& -2<K_0\xi,\eta>\\
c^2\rho_3(B(J_0\xi),B(\eta))+c^2\rho_3(B(\xi),B(J_0\eta)) &=& 2<I_0\xi,\eta>\nonumber
\end{eqnarray}
Put $J_1=I_0, J_2=J_0, J_3=K_0$ in (\ref{ric2}) we derive (\ref{ric3}).
Combining the already proved (\ref{ric3}) with (\ref{ric1})
we get that $\rho_{\alpha}$ is of type (1,1) with respect
to $J_{\alpha}$. Hence, the QKT structure is of instanton type.
The rest of b) follows by Theorem 6.1 and Remark 6.3 in \cite{PS1}.

We use (\ref{2.7}) again to prove c). We  have
$0=tr F_i(j_{*u}A^*)=-{c^2\over 2}\sum_{k=1}^{4n}(A,\Omega(B(e_k),B(J_0e_k)))$
which is equivalent to $ \sum_{k=1}^{4n}\rho_{\alpha}(e_k,J_{\beta}e_k)=0$
by Proposition~\ref{nov}. Further,
$0=tr F_i(j_{*u}B(\xi))=-{1\over 2}\sum_{k=1}^{4n}<\Theta(B(e_k),B(J_0e_k)),\xi>=t(u(\xi))$.

The proof of d) is a direct
consequence of the last equality in (\ref{2.7}) and (\ref{t1}). \hfill
{\bf Q.E.D.}

\

\medskip

\noindent
{\bf Remark 2.} We may consider the twistor space of an almost
quaternionic hermitian manifold and construct the almost complex
structure $I_1$ using horizontal spaces of a quaternionic connection
with skew-symmetric torsion. It follows that $I_1$ is integrable if
and only if the torsion is (1,2)+(2,1)-form, i.e. it is a QKT manifold.

\

\medskip

\noindent
{\bf Examples.} The twistor space $(Z,h_c,I_i), i=1,2$ of balanced HKT structures on the
nilpotent Lie groups constructed in \cite{DF} is semi-K\"ahler for $I_2$ and hermitian semi-K\"ahler
(balance) for $I_1$.

\section{Geometry of HKT manifold and twistor construction}

We recall some notations. The Lee form $\theta$ of a $2n$-dimensional
 almost Hermitian manifold $(M,g,J)$ with K\"ahler form
 $\Phi=g(\bullet,J\bullet )$ is defined by $\theta =-\delta \Phi\circ
 J$. On a HKT manifold there are three Lee forms corresponding to
 $J_{\alpha}, \alpha =1,2,3$, which are all equal. The common Lee form
 $\theta$ is called the Lee form of the HKT structure. It turns out
 that the Lee form of a HKT manifold is equal to the torsion 1-form,
 $\theta=t$ \cite{IP,I1}.

On a HKT manifold all the Ricci forms vanish and the exterior
differential of the Lee form is of type (1,1) with respect to each
$J_{\alpha}$ \cite{AI,IP}. The latter property can be easily seen
by comparing the curvatures of the Bismut and Chern connection
taken with respect to any hermitian structure $J_{\alpha}$ on the
corresponding canonical bundle. They differ by
$d(J_{\alpha}\theta)$, where
$J_{\alpha}\theta(X)=-\theta(J_{\alpha}X)$.  The curvature of the
Chern connection is of type (1,1), the curvature of Bismut
connection vanishes and the (2,0)+(0,2)-parts of $d\theta$ and
$d(J_{\alpha}\theta)$ coincide. Thus, Theorem~\ref{20} implies
that on a HKT manifold the Riemannian Ricci form $\rho^g_{\alpha}$
is of type (1,1) with respect to the complex structure
$J_{\alpha}$ and therefore the $*$-Ricci tensors are symmetric. We
have proved
\begin{th}\label{sim}
The $*$-Ricci tensors on a HKT manifold are symmetric.
\end{th}
\begin{pro}\label{prhkt}
Let $(Z,h_c,I_i, i=1,2)$ be a twistor space of a 4n dimensional
$(n>1)$ HKT manifold $(M,g,(J_{\alpha}))$ and
$X=u(\xi),Y=u(\eta)\in T_{\pi(u)}M$. The Ricci tensor $\rho$ and
the $*$-Ricci tensor $\rho^*_i$ for $(Z,h_c,I_i), i=1,2$ are given
by
\begin{eqnarray}\label{ri}\nonumber
& &Ric(j_{*u}A^*,j_{*u}B^*)=
\rho_{I_i}^{h_c}(j_{*u}A^*,j_{*u}B^*)=\frac{1}{nc^2}h_c(j_{*u}A^*,j_{*u}B^*),
\quad i=1,2;\\ \nonumber & &
Ric(j_{*u}B(\xi),j_{*u}A^*)=\rho_{I_i}^{h_c}(j_{*u}B(\xi),j_{*u}A^*)=0,\quad
i=1,2;\\ \nonumber & & Ric(j_{*u}B(\xi),j_{*u}B(\eta))= Ric^g(X,Y)
\\ \nonumber & &
\rho_{I_i}^{h_c}(j_{*u}B(\xi),j_{*u}B(\eta))=\rho_J^*(X,Y), \quad
J=j(u), \quad i=1,2. \nonumber
\end{eqnarray}
In particular, the $*$-Ricci tensor $\rho^*_i$ for $(Z,h_c,I_i)$
 is symmetric and $I_i$-invariant, $i=1,2$.

If the HKT space is Einstein (resp. $*$-Einstein with respect to
each $J_{\alpha}$) with positive scalar curvature $Scal^g$ (resp.
$Scal^g_Q$) then there exists an Einstein hermitian structure
$(Z,h_c,I_1), c^2=\frac{4}{Scal^g}$ (resp. $*$-Einstein almost
hermitian structure $(Z,h_c,I_2), c^2=\frac{4}{Scal^g_Q}$).
\end{pro}
{\bf Proof.}  Take the trace into (\ref{2.9}) and compare the
result with (\ref{r5}), (\ref{rn2}) to get the formulas in the
theorem. The last one and Theorem~\ref{sim} imply that the
$*$-Ricci tensors on the twistor space are symmetric. The formula
for the constant $c^2$ is a consequence of the fact that the
$*$-Einstein curvature is exactly equal to the quaternionic
curvature by Proposition~\ref{orp2}. \hfill {\bf Q.E.D.}\\

\noindent {\bf Remark 3}.  In view of the above results, the
$*$-Einstein condition on a HKT manifold does not impose
restrictions on the (2,0)+(0,2)-part of the $*$-Ricci tensor.\\

On a HKT manifold the quaternionic scalar curvature $Scal_Q=0$. In
this case Theorem~\ref{th11} leads to the following

\begin{th}
Let $M$ be a compact 4n-dimensional $(n>1)$ HKT manifold. Then
\begin{itemize}
\item[a)] $\int_M(Scal^g-Scal^g_Q)\,dV\ge 0$ with the equality if and only if the
HKT structure is balanced;
\item[b)] $\int_M(Scal^g-2Scal^g_Q)\,dV\ge 0$ with the equality if and only if the
HKT structure is hyperK\"ahler;
\end{itemize}
In particular, any compact HKT manifold with flat metric is
hyperK\"ahler.
\end{th}

On a HKT manifold $dT$ is (2,2)-form with respect to each complex
structure $J_{\alpha}$ and therefore all terms in the second line
in (\ref{22}) are equal (see \cite{I1}). Thus, (\ref{22}) gives
\beq\label{bal}
Ric(X,Y)=(\nabla_X\theta)Y+\frac{1}{4}\sum_{i=1}^{4n}
dT(X,J_{\alpha}Y,e_i,J_{\alpha}e_i).\eeq

 As a consequence of (\ref{bal}) we get that on a 4n
dimensional ($n\ge2$) HKT manifold the Ricci tensor is symmetric
if and only if $d^{\nabla}\theta(X,Y)=0.$ In particular, on a
balanced 4n dimensional ($n\ge2$) HKT manifold the Ricci tensor is
$J_{\alpha}$-invariant and symmetric. Therefore, the torsion
3-form is coclosed, $\delta T=0$. We note that this is true in
more general situation on any balanced Hermitian manifold which
Bismut connection has holonomy contained in the special unitary
group $SU(n)$ \cite{IP}.

\bigskip {\bf Authors' address:}
Stefan Ivanov, Ivan Minchev\\ University of Sofia, Faculty of Mathematics and
Informatics, Department of Geometry, \\ 5 James Bourchier blvd,
1164 Sofia, BULGARIA.\\ E-mail: {\tt ivanovsp@fmi.uni-sofia.bg}
\end{document}